\documentclass[12pt]{article}
\usepackage{a4wide}

\usepackage{comment}
\usepackage{multirow}
\usepackage{latexsym,amssymb,amsmath,epsfig,amsfonts}
\usepackage{algorithm}
\usepackage[noend]{algpseudocode}
\usepackage[title]{appendix}
\usepackage{eucal}
\usepackage{multirow,multicol,array,bm}
\usepackage{mathtools}
\usepackage{tabularx, booktabs, makecell, caption}
\usepackage{siunitx}
\usepackage[colorlinks=true,allcolors=blue]{hyperref}
\usepackage{float}
\usepackage[utf8]{inputenc}
\usepackage[english]{babel}
\usepackage{comment}
\usepackage{bbm}
\usepackage{graphicx}
\usepackage{float}
\usepackage[skip=0.333\baselineskip]{subcaption}
\usepackage[font=normal,labelfont=bf]{caption}

\numberwithin{table}{section}
\numberwithin{equation}{section}
\usepackage{amsthm}

\usepackage{color} 

\newtheorem{theorem}{Theorem}

\newtheorem{proposition}[theorem]{Proposition}

\newtheorem{definition}{Definition}

\newcommand{\jiesen}[1]{{{\color{black}#1}}}
\newcommand{\js}[1]{{{\color{black}#1}}} 
\newcommand{\camiel}[1]{{{\color{black}#1}}}

\begin{document}
\title{Strategic timing of arrivals to a queueing system with scheduled customers}

\author{Wathsala Karunarathne\thanks{Teletraffic Research Centre, School of Computer and Mathematical Sciences, The University of Adelaide, SA 5000, Australia, wathsala.karunarathne@adelaide.edu.au.} \and
Camiel M.P. Koopmans\thanks{Mathematical Institute, Leiden University, 2333 CC Leiden, The Netherlands, c.m.p.koopmans@math.leidenuniv.nl.} \and
Jiesen Wang \thanks{Korteweg-de Vries Institute, University of Amsterdam, 1090 GE Amsterdam, The Netherlands, j.wang2@uva.nl.}}

\maketitle

\begin{abstract} 

This paper examines a single-server queueing system that serves both scheduled and strategic walk-in customers. The service discipline follows a first-come, first-served policy, with scheduled customers granted non-preemptive priority.
Each walk-in customer strategically chooses their arrival time to minimise their expected waiting time, taking into account the reservation schedule and the decisions of other walk-in customers. \jiesen{We derive the Nash equilibrium arrival distribution for walk-in customers. We also study the case where early arrivals are allowed and analyse its impact on equilibrium arrival patterns.
By analysing various appointment schedules, we assess their effects on equilibrium arrival behaviour, waiting times, and server idle time. Finally, we compare different scheduling policies and evaluate their impact on overall system performance.}

\end{abstract}

\section{Introduction} \label{section:introduction}

Appointment scheduling systems are fundamental to service operations across numerous sectors, from healthcare to personal services like hairdressers. While these systems traditionally manage scheduled appointments, many service providers also accommodate walk-in customers, creating a complex dynamic between scheduled and opportunistic service seeking behaviour of walk-in customers. This dual access model presents unique challenges for both service providers and customers, particularly when walk-in customers behave strategically.

The complexity of dual access systems arises from several interrelated factors. First, service providers must balance the predictability of scheduled appointments with the flexibility demand of walk-in customers. Second, the system performance is significantly influenced by the strategic behaviour of walk-in customers, who make arrival decisions based on their understanding of the scheduling system and expected wait times. Third, the provider faces conflicting operational objectives; minimising server idle time while also minimising customer waiting times, and maintaining the service quality for all the customers while securing the system stability.

Consider a common scenario in personal services: a customer intending a visit to a hairdresser without an appointment. Such customers typically make arrival decisions based on available information, such as scheduled appointment slots and  historical wait times, which are affected by other customers' arrival decisions. This decision making process becomes naturally strategic, as customers aim to optimise their objectives considering other's decisions as well. Simultaneously, service providers place their own set of objectives, including maximising revenue and minimising server idle time.

The strategic behaviour of walk-in customers introduces a game theoretic element to the usual appointment scheduling problem. Each walk-in customer's optimal decision is not only affected by her own arrival time, but also by arrival times of other customers; both scheduled and walk-in. This creates a complicated feedback loop where individual decisions collectively form system performance, which in return influences customer decisions. Understanding and modelling this strategic behaviour is crucial for developing effective scheduling systems.

There are many advantages of optimal scheduling in dual access systems. In healthcare settings, effective management of scheduled appointments and walk-in patients not only has an impact on customer outcomes through reduced wait times and better access to care, but also improves resource utilisation. Commercial services benefit from the flexibility to serve both customer types while ensuring the systems' stability and customer satisfaction. Public organisations such as government offices and community facilities can better manage peak demand periods by balancing scheduled customers, leading to more efficient service delivery.

Although there is a significant body of research in appointment scheduling, the literature has primarily focused on either pure appointment systems or simplified walk-in models. The seminal work by Bailey and Welch~\cite{Welch52} followed by subsequent research established fundamental principles for appointment scheduling, but these traditional models lack the ability to capture the complexity of modern dual access systems, especially when strategic walk-in customers are involved. However, recent advancements in queueing theory and game theory open up the potential to understand strategic behaviour in service systems. The integration of strategic behaviour analysis with traditional scheduling methods tackles theoretical challenges while providing practical approaches for optimal service system design.

In this paper, we characterize the Nash equilibrium behaviour of walk-in customers given a schedule of pre-arranged customer arrivals. \js{We focus on symmetric Nash equilibria, in which all walk-in customers adopt the same strategy.} We then investigate how the schedule of pre-arranged arrivals affects the equilibrium decisions of walk-in customers and, in turn, overall system performance.
The main challenges are as follows:
\begin{itemize}
    \item The interaction between scheduled arrivals and walk-in customers induces a non-trivial waiting structure. Due to the absolute priority of scheduled customers, a walk-in customer may be postponed whenever a scheduled customer arrives before their service is completed. Crucially, this effect is stochastic and depends on two sources of randomness: the time it takes to serve all the customers ahead and the service time of the future scheduled customers themselves. Depening on this, a customer who initially faces a given position in the system may experience zero, one, or multiple additional delays.
    \item The structure of strategic arrival decisions is clear in the classical setting without scheduled customers, where arriving at time zero is a natural and optimal choice due to the absence of priority interference. However, in the presence of a scheduled customer at time zero who is given absolute priority, arriving at time zero may no longer guarantee a favorable position in the service order. Moreover, if early arrival before service commencement is allowed, customers do not necessarily have an incentive to arrive earlier. Additionally, walk-in customers typically do not wish to arrive immediately after the scheduled customer. Our numerical results also suggest that this decision is parameter-dependent.
\end{itemize}

The paper is organised as follows. Section \ref{section:literature} reviews the existing literature on appointment systems with walk-in customers and strategic behaviour, highlighting the gap in research on scheduling systems that account for strategic walk-in customers. \jiesen{Section \ref{section:model} presents the model in detail and introduces the quantities required for the computation of the Nash equilibrium arrival distribution.} Section \ref{section:Nash} derives the Nash equilibrium arrival distribution for walk-in customers in systems that either allow or restrict early arrivals. \js{Section 
\ref{section:algex} presents the algorithm and provides illustrative examples to 
demonstrate its implementation.} \jiesen{Section \ref{section:schedule} analyses the effects of different appointment schedules on system performance, including equilibrium arrivals, waiting times, and server idle time.} Finally, Section \ref{section:conclusion} concludes with a summary of the model and outlines potential directions for future research.

\section{Related Literature }\label{section:literature}

\jiesen{The scheduling design} of dual access service systems with mixed customer arrivals; both scheduled and walk-in customers represent a significant challenge in operations management, particularly when walk-in customers make their decisions strategically. The complexity of analysing such systems arises from two key aspects; deriving the arrival distribution of strategic walk-in customers \jiesen{and determining optimal appointment schedules for scheduled customers given the arrival behaviour of walk-in customers.}

Our work addresses a decision-making problem involving walk-in customers in a system that also accommodates scheduled customers. Scheduled customers are given priority for service, but service cannot be preempted once it begins. Walk-in customers behave strategically, choosing their arrival times while considering the decisions of other walk-in customers as well as the existing schedule. As a result, the arrival distribution of walk-in customers forms a Nash equilibrium. Our research connects to two key areas: systems that manage both scheduled and walk-in customers, and the study of strategic behaviour in queueing systems.

\subsection{Appointment systems with walk-in customers}

Bailey and Welch~\cite{Welch52} pioneered the mathematical foundation of appointment scheduling by analysing the fundamental trade-off between patient waiting time and physician idle time. Their study of outpatient departments demonstrated that consultation time variability plays a critical role in queue formation and system performance. Subsequently, many researchers have studied appointment scheduling, mainly in healthcare systems~\cite{Ho92,Cayirli03,Robinson03,Kaandorp07}. Research on appointment scheduling has primarily focused on systems with only scheduled customers, considering both static scheduling~\cite{Pegden90,Stein94,Hassin08,Kuiper15} and dynamic scheduling~\cite{Wang93,Hahn14,Gilbertson16,Mahes21,Hassin24}. Moreover, real-world complications in appointment systems have received significant attention. Studies have investigated the impact of no-shows~\cite{Robinson10}, cancellations~\cite{Hahn14}, and early and late arrivals~\cite{Koeleman12}, leading to the development of more robust scheduling policies.

The above-mentioned studies focus on systems that accept customers only by prior appointment. However, many systems also accept walk-in customers to increase revenue and expand their customer base, while rejecting walk-in customers may result in a loss of goodwill. In their review of the literature on healthcare scheduling, Cayirli and Veral~\cite{Cayirli03} noted that walk-in patients are often neglected in scheduling studies, despite their presence in practice. Gupta, Zoreda, and Kramer~\cite{Gupta71} analysed the challenges of staff scheduling in hospitals in the presence of walk-in patients. Fiems, Koole, and Nain~\cite{Fiems07} studied the effect of emergency arrivals on scheduled patients' waiting times. Robinson and Chen~\cite{Robinson10} compared two appointment scheduling systems: a traditional system in which routine patients book appointments well in advance but have relatively high no-show rates, and an open-access system serving same-day patients. They found that open-access scheduling generally outperforms traditional scheduling because the variability caused by no-shows can be more disruptive than the variability in same-day patient arrivals. Alexandrov and Lariviere~\cite{Alexandrov12} examined whether restaurants should offer reservations under uncertain demand. They found that reservations can increase sales during slow periods by guaranteeing seats to potential diners who might otherwise stay home.

Although walk-in arrivals can reduce system efficiency, rejecting them may not be desirable in many practical settings. Despite the prevalence of this problem, the relevant literature remains relatively limited. Rahmani and Heydari~\cite{Rahmani14} proposed a two-step proactive-reactive scheduling approach for a two-machine flow shop system under uncertainty. Their approach balances the trade-off between accounting for uncertain processing times in the initial schedule and adapting the schedule to unexpected job arrivals during execution. Hassin and Wang~\cite{Hassin24} considered a healthcare system with both urgent and routine patients and studied how many appointment slots should be allocated to routine patients. They proposed a two-cutoff heuristic that performs close to optimally when routine patient arrivals approach system capacity. Their results show that traditional cutoff strategies can become inefficient when the mean number of routine patients approaches daily capacity, whereas their proposed heuristic provides a better approximation. Wang, Liu, and Wan~\cite{Wang20} developed a static scheduling model in the presence of walk-in customers. Their objective was to determine the number of customers to schedule in each equally spaced time slot so as to minimize the total expected waiting cost of customers, server idle cost, and overtime cost. Karunarathne~\cite{Karunarathne23} extended this line of research by determining the arrival times of scheduled customers in the presence of walk-in customers. The objective was to minimize a linear combination of customers' total expected waiting cost and the server's total expected idle cost, with scheduled customers receiving non-preemptive priority.

\subsection{Strategic behaviour of walk-in customers}

The strategic behaviour of customers in queues, particularly in settings involving strategic interactions and the distinction between private and social optimization, has received considerable attention in Operations Research, Management Science, Economics, and Computer Science. Since the foundational work of Naor~\cite{Naor69}, a substantial body of literature has developed on these topics. Comprehensive reviews can be found in Hassin and Haviv~\cite{Hassin03,Hassin16}.

The strategic decision of interest in this paper is when to arrive. Glazer and Hassin~\cite{GH83} examined a single-server queue in which customers independently choose their arrival times. The system operates from time $0$ to time $T$ on a first-come, first-served (FCFS) basis, with exponentially distributed service times. Customers seek to minimize their expected waiting times, the total number of arrivals follows a Poisson distribution, and customers may arrive before the system opens at time $0$.

Haviv~\cite{Haviv13} considered a single-server system with independent, exponentially distributed service times in which walk-in customers seek to minimize the sum of tardiness and waiting costs. He used numerical procedures to characterize equilibrium arrival patterns and fluid models to derive explicit equilibria. Extending this line of research, Sherzer and Kerner~\cite{Sherzer17} incorporated earliness costs into the model. Their results showed that the symmetric equilibrium is unique and mixed when customers have complete freedom to choose their arrival times. Haviv and Ravner~\cite{Haviv15} studied a $M/M/m/c$ queue with strategic walk-in customers who choose their arrival times. They analysed loss probabilities both in a game-theoretic setting, in which customers act selfishly, and under social optimization. They showed that the equilibrium arrival distribution is not generally uniform when the customer population is insufficiently large. In their survey of queueing systems with strategic timing of arrivals, Haviv and Ravner~\cite{Haviv21} reviewed mathematical modelling approaches and identified open questions concerning customers' arrival-time decisions. Their survey particularly discussed the ratio between socially optimal and equilibrium outcomes, commonly known as the price of anarchy, which measures the inefficiency caused by strategic decision-making. For a comprehensive review of this literature, see Chapter~4 of Hassin and Haviv~\cite{Hassin16}.

Cil and Lariviere~\cite{CL13} studied a service system in which customers can either reserve seats in advance or walk in, choosing between receiving service and staying home. Tunçalp et al.~\cite{Tuncalp24} considered a clinic in which strategic patients choose whether to schedule appointments or walk in, while the clinic maximizes revenue by allocating appointment slots between the two groups. Walk-in patients may be denied service. Baron, Chen, and Li~\cite{BCL23} examined a system in which customers can pre-order online before visiting a store or walk in directly, with walk-in customers deciding whether to wait or balk. Liu, van Jaarsveld, Wang, and Xiao~\cite{LJJWX23} investigated an outpatient system serving both routine and urgent patients, who can book appointments, walk in, or balk. They examined the effects of service-capacity allocation, disclosure of appointment delays, and restrictions on walk-in availability.

To the best of the author's knowledge, no prior research has explored the topic of strategic arrivals in the presence of scheduled customers. This paper addresses this gap by developing a comprehensive model that examines the arrival distribution of walk-in customers and the influence of the appointment schedule on this distribution. Specifically, we derive the equilibrium arrival distribution of strategic walk-in customers for a given schedule, capturing how these customers optimise their arrival timing decisions. Then we analyse the impact of various scheduling schemes on the arrival distribution, and consequently, on the expected waiting time for each customer and the expected idle time of the system.


\section {Model and preliminaries} \label{section:model}

The system begins operation at time $0$ and closes at time $T$. 
Each day, the system serves $M$ scheduled customers and a Poisson distributed number of walk-in customers, with a mean of $\lambda$. \camiel{The service time of both scheduled customers and walk-in customers is i.i.d. and exponentially distributed with rate $\mu$.} The Poisson distribution for walk-in customers has a number of favorable properties that simplify analysis (see \cite{Haviv12}). \jiesen{A key simplification is that arrivals in disjoint time intervals are independent. This property facilitates the separate analysis of different time segments.
More specifically, by the thinning property, the number of customers joining within any given time interval remains Poisson distributed. Also, when focusing on a tagged customer, the remaining walk-in customers retain the Poisson structure with rate \(\lambda\).}   
To make the appointment schedule, we need to decide $M$ distinct arrival times for all $M$ customers on $[0,T]$. In the framework we consider, it is not allowed to use randomised schedules.
The service operates under a first-come, first-served (FCFS) discipline, with scheduled customers receiving non-preemptive priority over walk-in customers. Non-preemptive priority means that once a service begins, it cannot be interrupted by the arrival of a higher-priority customer. If a scheduled customer arrives while the system is busy, she will wait until the current service is completed and all previously scheduled customers have been served. \jiesen{When multiple walk-in customers arrive at the same time, the order between these customers is decided uniformly at random.} The system closes at time $T$, after which no walk-in customers are accepted. However, all customers who arrived before $T$ will still be served.

Walk-in customers can choose their arrival times to minimise their expected waiting time. We will consider both the restriction of possible arrival times to the domain $[0,T]$ and to the domain $(-\infty,T]$. We refer to arrivals before time 0 as early arrivals.   The decision of arrival time is influenced by the appointment schedule and the choices of other walk-in customers. We focus on a symmetric strategy, where all walk-in customers adopt the same approach. \jiesen{This strategy is represented as a probability distribution over the interval $[0,T]$ or $[-\infty,T]$, depending on whether early arrivals are allowed.}

The appointment schedule for $M$ customers is represented as a vector $\mathcal{T}_s\in \mathbb{R}_{\geq 0}^{M}$, where $\mathcal{T}_s(i)$ corresponds to a scheduled arrival time of customer $i$, for $i= 1,\dots, M$. We assume that there will be at most one customer scheduled at one time point, that is, $\mathcal{T}_s(i) < \mathcal{T}_s(i+1)$ for $i = 1,\ldots, M-1$. Let
$\mathcal{T}\in \mathbb{R}_{\geq 0}^{M+2}$ be an augmented vector based on $\mathcal{T}_s$ where the boundaries of the schedule are given by $\mathcal{T}(0)=0, \mathcal{T}(M+1)=T$, and $\mathcal{T}(i) = \mathcal{T}_s(i) $ for $i= 1,\dots, M$.
It is possible to schedule a customer at time $0$ or $T$, in which case $\mathcal{T}(1) = \mathcal{T}(0)$ or $\mathcal{T}(M) = \mathcal{T}(M+1)$.

\js{Let $F(t)$ denote the arrival-time distribution of walk-in customers, with density $f(t)$. We allow $F(t)$ to have finitely many atoms. Our objective is to characterize the Nash equilibrium arrival distribution. Given an arrival distribution $F$, let $N(t)$ denote the number of customers in the system at time $t$, and define
\[
    P_n(t)=\mathbb{P}(N(t)=n).
\]
At each scheduled arrival time $\mathcal{T}(k)$, the system gains one customer. Accordingly, the state probabilities satisfy
\begin{equation}\label{eq:PschedArr}
    P_{n+1}(\mathcal{T}(k))
    = \lim_{t\uparrow \mathcal{T}(k)} P_n(t),
    \qquad n\geq 0.
\end{equation}
Between scheduled arrival times, the state probabilities evolve according to the forward equations
\begin{align*}
    \frac{\partial}{\partial t} P_0(t)
    &= \mu P_1(t)-\lambda f(t)P_0(t),
    \\
    \frac{\partial}{\partial t} P_n(t)
    &= \lambda f(t)P_{n-1}(t)+\mu P_{n+1}(t)
    -\bigl(\lambda f(t)+\mu\bigr)P_n(t),
    \qquad n=1,2,\ldots.
\end{align*}

In addition, let $w_{\mathcal{T},k}(n,t)$ denote the expected waiting time of a walk-in customer arriving at time $t$, conditional on the system containing $n$ customers, under the schedule $\mathcal{T}$. This conditional expected waiting time also depends on $F$, since the arrival distribution of other walk-in customers affects the subsequent evolution of the system.

The expected waiting time of a walk-in customer arriving at time $t$ can therefore be expressed as
\begin{equation} \label{eq:Ewt}
    E_w(t):= \sum_{n=0}^\infty P_n(t)w_{\mathcal{T},k}(n,t).
\end{equation}
Characterizing $w_{\mathcal{T},k}(n,t)$ is particularly challenging because scheduled customers have absolute priority and may postpone a walk-in customer's service even after she has arrived. Thus, the waiting time depends not only on the number of customers already in the system, $n$, but also on the timing of subsequent scheduled arrivals and the possibility of repeated postponements. The index $k$ specifies the scheduled arrivals that need to be considered in calculating the expected waiting time. Strictly speaking, $k$ can be determined from $\mathcal{T}$ and $t$; however, we retain $k$ in the notation to make the interpretation of the waiting-time function more transparent.
}

{\bf Calculation of \js{$w_{\mathcal{T},k}(n,t)$}.} \js{The expected waiting time can be calculated recursively. At the first step, there are two possible cases. First, the walk-in customer is served before the next scheduled customer arrives at $\mathcal{T}(k+1)$. In this case, no additional scheduled customer needs to be considered. Second, the walk-in customer is not served before $\mathcal{T}(k+1)$. In this case, the customer must wait for the additional scheduled customer arriving at $\mathcal{T}(k+1)$. Moreover, we must further consider whether the walk-in customer can be served before the following scheduled customer arrives at $\mathcal{T}(k+2)$. The analysis then proceeds in the same manner: at each subsequent step, we consider whether the walk-in customer is served before the next scheduled arrival. This gives rise naturally to the recursive formulation.

The number of customers served during an interval of length $x$ follows a Poisson distribution with parameter $\mu x$. We denote the probability of serving $i$ customers during this interval by
\[
    \operatorname{Pois}(x;i)
    := \frac{e^{-\mu x}(\mu x)^i}{i!}.
\]
Let $X$ follow an Erlang distribution with shape parameter $n$ and rate parameter $\mu$. Its cumulative distribution function and probability density function are
\[
    F_{\mathrm{Erl}}(x;n,\mu)
    = 1-\sum_{i=0}^{n-1}\operatorname{Pois}(x;i),
\]
and
\[
    f_{\mathrm{Erl}}(x;n,\mu)
    = \frac{\mu^n x^{n-1}e^{-\mu x}}{(n-1)!}
    = \mu\operatorname{Pois}(x;n-1),
\]
respectively.

\paragraph{Step 1: Detailed representation of the expected waiting time.}
For $t\geq 0$, consider a walk-in customer arriving at time $t$, when there are $n\geq 1$ customers ahead in the queue. Let
\[
    h=\mathcal{T}(k+1)-t
\]
denote the time until the next scheduled arrival. We condition on the number of customers served during $(t,\mathcal{T}(k+1))$.

If all $n$ customers are served before the next scheduled arrival, the walk-in customer begins service after an Erlang-distributed time. The first term below gives the expected waiting time contributed by this event. Otherwise, if exactly $i<n$ customers are served during the interval, the walk-in customer must wait until time $\mathcal{T}(k+1)$ and then wait for the remaining customers ahead in the queue. Since the number of customers served during the interval follows a Poisson distribution, we obtain
\begin{align}
    w_{\mathcal{T},k}(n,t)
    &= \int_0^{h} x f_{\mathrm{Erl}}(x;n,\mu)\,dx \nonumber\\
    &\quad + \sum_{i=0}^{n-1}
    \operatorname{Pois}(h;i)
    \left(
        h+w_{\mathcal{T},k+1}
        \bigl(n-i+1,\mathcal{T}(k+1)\bigr)
    \right).
    \label{eq:wkt}
\end{align}
The integral accounts for the case in which all $n$ customers are served before the next scheduled arrival, while the summation accounts for the cases in which fewer than $n$ customers are served.

\paragraph{Step 2: A tractable representation.}
To simplify \eqref{eq:wkt}, define
\[
    \widetilde{w}_{\mathcal{T},k}(n)
    := w_{\mathcal{T},k}\bigl(n,\mathcal{T}(k)\bigr)
\]
and let
\begin{align*}
    \mathcal{E}(x;n,\mu)
    &:= \int_0^x y f_{\mathrm{Erl}}(y;n,\mu)\,dy\\
    &= \frac{1}{(n-1)!}\int_0^x
        (\mu y)^n e^{-\mu y}\,dy\\
    &= \frac{\gamma(n+1,\mu x)}{\mu(n-1)!},
\end{align*}
where $\gamma$ denotes the lower incomplete gamma function. Thus, $\mathcal{E}(x;n,\mu)$ is the expected waiting-time contribution from service completion before time $x$. Substituting these definitions into \eqref{eq:wkt} yields
\begin{align}
    w_{\mathcal{T},k}(n,t)
    &= \mathcal{E}\bigl(\mathcal{T}(k+1)-t;n,\mu\bigr)
    \nonumber\\
    &\quad + \sum_{i=0}^{n-1}
    \operatorname{Pois}\bigl(\mathcal{T}(k+1)-t;i\bigr)
    \left(
        \mathcal{T}(k+1)-t
        +\widetilde{w}_{\mathcal{T},k+1}(n-i+1)
    \right).
    \label{eq:wkt_compact}
\end{align}
For completeness, $w_{\mathcal{T},k}(0,t)=0$, since no customers are ahead of the arriving walk-in customer when $n=0$. Moreover, $w_{\mathcal{T},k}(n,t)$ is continuous in $t$ between scheduled arrival times. At a scheduled arrival time, the waiting-time function satisfies
\[
    \lim_{t\uparrow\mathcal{T}(k+1)}
    w_{\mathcal{T},k}(n,t)
    =
    w_{\mathcal{T},k+1}
    \bigl(n+1,\mathcal{T}(k+1)\bigr).
\]
If we let $t_k: = \mathcal{T}(k+1) - \mathcal{T}(k)$ be the inter-arrival times for scheduled customers.
The value of $\Tilde{w}_{\mathcal{T},k}(n)$ can be determined through Equation~(\ref{eq:wkt}) as
\begin{align*}
     & \Tilde{w}_{\mathcal{T},k}(n)= \mathcal{E}(t_k;n,\mu) + \sum_{i=0}^{n-1} 
     \text{Pois}(t_k; i)  \left(t_k+\Tilde{w}_{\mathcal{T},k+1}(n-i+1)\right), \\
     & \Tilde{w}_{\mathcal{T},M}(n) = \frac{n}{\mu} \,. 
\end{align*}

The above quantities jointly determine the expected waiting time
process $E_w(t)$: the arrival distribution $F$ determines the state
probabilities $P_n(t)$, while $P_n(t)$ and the conditional expected
waiting times $w_{\mathcal{T},k}(n,t)$ together determine the
unconditional expected waiting time. In Section~4, we use this
relationship to characterize the Nash equilibrium arrival
distribution by imposing the condition that no walk-in customer
can reduce her expected waiting time by unilaterally changing her
arrival time.
}

\section{Nash equilibrium} \label{section:Nash}


\js{Let $f_e$ and $F_e$ denote the probability density function and cumulative distribution function, respectively, of the Nash equilibrium arrival distribution. We allow $F_e$ to have both a continuous component and finitely many atoms, each representing an arrival time with positive probability mass. However, we will show later that the equilibrium distribution can have at most one atom. We will also explain how to identify this atom when it exists.

In a Nash equilibrium, all arrival times in the support of $F_e$ must yield the same expected waiting time when all walk-in customers follow $F_e$. Otherwise, a customer could unilaterally deviate to another arrival time and obtain a lower expected waiting time, contradicting the equilibrium condition. We denote this common expected waiting time by $E_w$. The equilibrium condition can therefore be expressed as follows.}

\begin{definition}
An arrival distribution $F_e$ is an \emph{equilibrium arrival distribution} \js{if} there exists a value $E_w \in \mathbb{R}_{\ge 0}$ such that the induced expected waiting time $E_w(t)$ satisfies
\[
E_w(t) = E_w, \quad t \in \operatorname{supp}(F_e),
\]
and
\[
E_w(t) \ge E_w, \quad t \notin \operatorname{supp}(F_e).
\]
\end{definition}
The equilibrium arrival distribution defined above is a mixed-strategy Nash equilibrium, with customers randomizing their arrival times according to $F_e$. Since all arrival times in the support of $F_e$ yield the same expected waiting time, while any arrival time outside the support yields a no smaller expected waiting time, no customer can benefit from a unilateral deviation. Therefore, $F_e$ is a best response to itself and thus constitutes a Nash equilibrium.

\jiesen{

We assume that customers follow the equilibrium strategy induced by $F_e$. 
Let $t \in [\mathcal{T}(k), \mathcal{T}(k+1))$ be such that $f_e(t) > 0$, i.e., $t$ lies in the interior of the support of $F_e$. Then, by the equilibrium (indifference) condition, all arrival times in the support yield the same expected payoff. In particular, we have
\begin{equation}\label{eq:Ew}
E_w(t) = \sum_{n=0}^\infty P_n(t)\, w_{\mathcal{T},k}(n,t)
= E_w,
\end{equation}
where $E_w$ is a constant independent of $t$, and $P_n(t)$ denotes the queue-length distribution induced by the arrival process $F_e$.}

\js{
To characterize the equilibrium arrival distribution $F_e$, we exploit the indifference condition that $E_w(t)$ is constant on the support of $F_e$. Thus, at interior points of the support where $E_w(t)$ is differentiable, the equilibrium condition requires
\[
\frac{\partial E_w(t)}{\partial t}=0.
\]
This condition provides a necessary condition for characterizing the equilibrium density $f_e(t)$, together with the requirement that no arrival time outside the support yields a lower expected waiting time.
The main difficulty in calculating $\partial E_w(t)/\partial t$ stems from the dependence of $w_{\mathcal{T},k}(n,t)$ on $t$; see Equation~(\ref{eq:DeGP}) in Appendix~\ref{Derivatives}. We therefore first derive an expression for $\partial w_{\mathcal{T},k}(n,t)/\partial t$.

To study how $w_{\mathcal{T},k}(n,t)$ changes with the arrival time $t$,
we derive its time derivative in several steps. Recall that
$h=\mathcal{T}(k+1)-t$, so that $\partial h/\partial t=-1$.
We first note that
\[
    \frac{\partial \operatorname{Pois}(h;0)}{\partial h}
    =-\mu\operatorname{Pois}(h;0),
\]
and, for $i\geq 1$,
\[
    \frac{\partial \operatorname{Pois}(h;i)}{\partial h}
    =\mu\bigl(\operatorname{Pois}(h;i-1)
    -\operatorname{Pois}(h;i)\bigr).
\]
We first differentiate the integral term in the expression for
$w_{\mathcal{T},k}(n,t)$. By the Leibniz integral rule,
\[
    \frac{\partial}{\partial t}
    \int_0^h x f_{\mathrm{Erl}}(x;n,\mu)\,dx
    =-h f_{\mathrm{Erl}}(h;n,\mu).
\]
The negative sign arises because the remaining time $h$ decreases as the walk-in customer arrives later.

Next, consider the summation term in the expression for
$w_{\mathcal{T},k}(n,t)$. For convenience, define
\[
    A_i(t):=
    h+\widetilde{w}_{\mathcal{T},k+1}(n-i+1),
    \qquad i=0,\ldots,n-1.
\]
Here, $A_i(t)$ is the waiting time at the next scheduled arrival, including the time $h$ until that arrival, conditional on exactly $i$ service completions during the interval. Since $\widetilde{w}_{\mathcal{T},k+1}(n-i+1)$ does not depend on $t$, we have $\partial A_i(t)/\partial t=-1$. Applying the product rule and the Poisson derivative identities gives
\begin{align}
    \frac{\partial}{\partial t}
    \sum_{i=0}^{n-1}
    \operatorname{Pois}(h;i)A_i(t) 
    & =
    -\sum_{i=0}^{n-1}\operatorname{Pois}(h;i)
    -\mu\operatorname{Pois}(h;0)A_0(t) \nonumber\\
    &\qquad
    -\mu\sum_{i=1}^{n-1}
    \bigl(\operatorname{Pois}(h;i-1)
    -\operatorname{Pois}(h;i)\bigr)A_i(t).
    \label{eq:derwkt_sum}
\end{align}
The first term reflects the reduction in the time remaining until the next scheduled arrival. The other terms capture how the probabilities of the different numbers of service completions change as the walk-in arrival time varies.

Combining the derivatives of the integral and summation terms, we obtain
\begin{align}
    \frac{\partial w_{\mathcal{T},k}(n,t)}{\partial t}
    &=-h f_{\mathrm{Erl}}(h;n,\mu)
    -\sum_{i=0}^{n-1}\operatorname{Pois}(h;i) \nonumber\\
    &\quad
    -\mu\operatorname{Pois}(h;0)
    \bigl(h+\widetilde{w}_{\mathcal{T},k+1}(n+1)\bigr)
    \nonumber\\
    &\quad
    -\mu\sum_{i=1}^{n-1}
    \bigl(\operatorname{Pois}(h;i-1)
    -\operatorname{Pois}(h;i)\bigr)
    \bigl(h+\widetilde{w}_{\mathcal{T},k+1}(n-i+1)\bigr).
    \label{eq:derwkt}
\end{align}
This derivative describes how the expected waiting time conditional on the number of customers ahead changes with the walk-in arrival time. A later arrival shortens the time until the next scheduled arrival, but it also changes the distribution of the number of customers served before that arrival. Consequently, the derivative
need not be negative: its sign depends on the balance between these effects.

We can now differentiate the unconditional expected waiting time,
\[
    E_w(t)=\sum_{n=0}^{\infty}P_n(t)w_{\mathcal{T},k}(n,t).
\]
Using Equation \eqref{eq:DeGP} in Appendix \ref{Derivatives}, we obtain
\begin{align}
    \frac{\partial E_w(t)}{\partial t}
    &=\sum_{n=0}^{\infty}P_n(t)
    \frac{\partial w_{\mathcal{T},k}(n,t)}{\partial t}
    \nonumber\\
    &\quad+\lambda f(t)\sum_{n=0}^{\infty}P_n(t)
    \bigl(w_{\mathcal{T},k}(n+1,t)
    -w_{\mathcal{T},k}(n,t)\bigr) \nonumber\\
    &\quad-\mu\sum_{n=0}^{\infty}P_{n+1}(t)
    \bigl(w_{\mathcal{T},k}(n+1,t)
    -w_{\mathcal{T},k}(n,t)\bigr).
    \label{eq:PartialEw}
\end{align}
The first term captures the direct effect of changing the arrival time on the conditional expected waiting time. The second term captures the effect of additional walk-in arrivals on the queue length distribution, while the third captures the effect of service completions. Together, these terms determine how the expected waiting time changes with the walk-in arrival time.
Imposing the equilibrium condition $\partial E_w(t)/\partial t=0$ and solving for $f(t)$ yields the following expression for the equilibrium arrival density:}
\begin{equation} \label{eq:density}
    \begin{split}
        f_e(t) &=  \frac{\mu \sum_{n = 0}^{\infty} P_{n+1}(t) \, (w_{\mathcal{T},k}(n+1,t) - w_{\mathcal{T},k}(n,t)) -  \sum_{n = 0}^\infty P_n(t) \,  \frac{\partial}{\partial t}w_{\mathcal{T},k}(n,t)}{\lambda \sum_{n = 0}^{\infty} P_n(t) \, (w_{\mathcal{T},k}(n+1,t) - w_{\mathcal{T},k}(n,t))}       \,.
    \end{split}
\end{equation}
The equilibrium arrival strategy is derived under the assumption that $F_e$ is continuous on $[\mathcal{T}(k), \mathcal{T}(k+1))$. We then verify ex post that the resulting equilibrium solution is consistent with this assumption, thereby confirming its validity.


\subsection{Not including early arrivals}\label{sect:notEarly}

\js{We first consider the case in which there are no arrivals before time $0$. This model setting admits two interpretations. First, early arrivals may be prohibited, so customers can join the queue only after the system begins operating at time $0$. Alternatively, early arrivals may be permitted, but arriving before time $0$ does not secure a customer's position in the queue. In the latter case, customers have no incentive to arrive before time $0$, since doing so provides no advantage over arriving when the system starts operating.}

The Nash equilibrium structure may differ depending on whether the first scheduled customer arrives at time \(0\). When there is no scheduled customer at time \(0\), walk-in customers have an incentive to arrive at time \(0\). However, when a scheduled customer arrives at time \(0\), this incentive may disappear due to the absolute priority given to scheduled customers.

\jiesen{We first claim that \(F_e\) does not have any atom at any time \(t \in (0,T]\). Indeed, if an atom existed at some positive time \(t\), then it would be strictly better for a customer to arrive slightly earlier than \(t\), thereby securing a better position in the queue. Hence, no atom can occur on \((0,T]\). At time \(0\), however, the situation is different. When \(\mathcal{T}(1)>0\), there is necessarily an atom at \(0\), since the first customer to arrive before the first scheduled customer would have waiting time 0. 
If there is no atom at time \(0\), arriving at time \(0\) with probability one would be a profitable deviation. When \(\mathcal{T}(1)=0\), there can be an atom at 0. Although the scheduled customer has absolute priority at time \(0\), arriving at time \(0\) can still improve a random customer’s position relative to other random arrivals.}

Note that when there is an atom \(p_0\) at time \(0\), the number of walk-in arrivals at \(t=0\) follows a Poisson distribution with parameter \(p_0\lambda\). \jiesen{However, the computation of the expected waiting time for a tagged customer arriving at \(t=0\) is more involved than in the standard model without scheduled customers. In particular, one must account for the possibility that the tagged customer does not enter service before the next scheduled customer arrives. In this case, the scheduled customer immediately joins the system with absolute priority, effectively adding one additional customer ahead of the tagged customer. This mechanism may repeat at subsequent scheduled arrival times. Consequently, when \(\mathcal{T}(1)>0\), the tagged customer may be delayed by at most \(M\) additional customers, whereas in the boundary case \(\mathcal{T}(1)=0\), where the first scheduled customer is already present upon the tagged customer’s arrival, the maximum delay is \(M-1\).

The expected waiting time of arriving at $0$ is therefore given by}
\begin{equation}\label{eq:Ew0}
   E_w(0) := \sum_{n=0}^\infty \left(\sum_{i = n}^\infty \frac{(p_0\lambda)^i}{i!}e^{-p_0\lambda}\frac{1}{i+1}\right)  w_{\mathcal{T},\mathbbm{1}_{\mathcal{T}(1)=0} }(n+\mathbbm{1}_{\mathcal{T}(1)=0},0) \,.
\end{equation}
\jiesen{
To interpret the formula, we consider the events that lead to \(n\) walk-in customers being in front of the tagged customer in the queue at time \(0\). At least \(i\ge n\) walk-in customers must have arrived simultaneously with her. Given this, with probability \(\frac{1}{i+1}\), there are exactly \(n\) walk-in customers ahead of her.
If a scheduled customer arrives at time \(0\), then \(\mathbbm{1}_{\{\mathcal{T}(1)=0\}}=1\). In this case, the scheduled customer occupies the first position in the queue, so the walk-in customers are indexed starting from \(1\), and the tagged customer has \(n+1\) customers ahead of her. Otherwise, \(\mathbbm{1}_{\{\mathcal{T}(1)=0\}}=0\), the walk-in customers are indexed starting from \(0\), and the tagged customer has \(n\) customers ahead of her.
}

\jiesen{We next claim that there is a gap immediately after each scheduled arrival time
$\mathcal{T}(k)$ when $\mathcal{T}(k)>0$: arriving at time
$t=\mathcal{T}(k)$ or $t=\mathcal{T}(k)+\delta$ is strictly worse than
arriving at time $t=\mathcal{T}(k)-\delta$ for sufficiently small
$\delta>0$. The system state is essentially identical in all three cases,
except that at time $\mathcal{T}(k)$ an additional scheduled customer with
absolute priority joins the system. An arriving walk-in customer at time $\mathcal{T}(k)-\delta$ can potentially start service before the scheduled customer, in contradiction with walk-in customers at time $\mathcal{T}(k)$ or $\mathcal{T}(k)+\delta$.

Note that this argument does not apply when $\mathcal{T}(1)=0$. In this case, a scheduled customer with absolute priority arrives at time zero. Arriving at $t=0$ coincides with this priority arrival, which places an arriving walk-in customer behind and in competition with other walk-in customers arriving at the same instant. This will be illustrated in the numerical example in Section \ref{section:algex}.}

\subsection{Including early arrivals}\label{sect:Early}

In this section, we enable the walk-in customers to arrive before time 0. This option leads to walk-in customers having the option to arrive before the system starts operating, with the advantage that a better position in the queue can be acquired. 
\jiesen{Similar to the previous case, if the first scheduled customer is at time \(0\), walk-in customers may not enter early in equilibrium, even when early entry is allowed.}

Due to the nature of the poisson distribution, the accumulated probabilities of customers in the system at time $t\leq 0$ is also Poisson distributed and given by
$$ P_n(t)=\frac{(\lambda F_e(t))^n e^{-\lambda F_e(t)}}{n!}.$$
\jiesen{ For an early customer arriving at time $t<0 $, with $n\in \mathbb{N}_0$ customers in front of her, the waiting time 
 is given by
\begin{equation*}
    \begin{split} w_{\mathcal{T},\mathbbm{1}_{\mathcal{T}(1)=0}}(n+\mathbbm{1}_{\mathcal{T}(1)=0},0)-t \, . 
    \end{split}
\end{equation*}
The role of $\mathbbm{1}_{\mathcal{T}(1)=0}$ is to add the scheduled customer at time 0 if there is one.

Suppose \(t\) belongs to the support of the equilibrium arrival distribution, i.e., \(f_e(t)>0\). Then, by the equilibrium condition, all arrival times in the support must yield the same expected waiting time. Let
\begin{equation*}
\begin{split}
E_w(t)
&:=\sum_{n=0}^{\infty}P_n(t)\,
w_{\mathcal{T},\mathbbm{1}_{\mathcal{T}(1)=0}}
\left(n+\mathbbm{1}_{\mathcal{T}(1)=0},t\right)\\
&=
-t+\sum_{n=0}^{\infty}P_n(t)\,
w_{\mathcal{T},\mathbbm{1}_{\mathcal{T}(1)=0}}
\left(n+\mathbbm{1}_{\mathcal{T}(1)=0},0\right)
=E_w,
\end{split}
\end{equation*}
where \(P_n(t)\) is computed under the equilibrium arrival density \(f_e\). Since \(E_w(t)\) is constant over the support of \(f_e\), differentiating both sides with respect to \(t\) and applying Equation ~\eqref{eq:deNPEarly1} gives
\begin{align*}
0
&=\frac{\partial E_w(t)}{\partial t}\\
&=
-1
+\frac{\partial}{\partial t}
\sum_{n=0}^{\infty}
P_n(t)
w_{\mathcal{T},\mathbbm{1}_{\mathcal{T}(1)=0}}
\left(n+\mathbbm{1}_{\mathcal{T}(1)=0},0\right)\\
&=
-1
+\lambda f_e(t)
\sum_{n=0}^{\infty}
\left(
w_{\mathcal{T},\mathbbm{1}_{\mathcal{T}(1)=0}}
\left(n+1+\mathbbm{1}_{\mathcal{T}(1)=0},0\right)
-
w_{\mathcal{T},\mathbbm{1}_{\mathcal{T}(1)=0}}
\left(n+\mathbbm{1}_{\mathcal{T}(1)=0},0\right)
\right)
P_n(t).
\end{align*}
Hence,
\begin{equation}\label{eq:earlyfe}
f_e(t)=
\frac{1}
{\lambda
\displaystyle\sum_{n=0}^{\infty}
\left(
w_{\mathcal{T},\mathbbm{1}_{\mathcal{T}(1)=0}}
\left(n+1+\mathbbm{1}_{\mathcal{T}(1)=0},0\right)
-
w_{\mathcal{T},\mathbbm{1}_{\mathcal{T}(1)=0}}
\left(n+\mathbbm{1}_{\mathcal{T}(1)=0},0\right)
\right)
P_n(t)}.
\end{equation}
Combining Equation~(\ref{eq:earlyfe}) for \(t<0\) with Equation~(\ref{eq:density}) for \(t\ge0\), we can obtain the equilibrium arrival density for this case.
}

\jiesen{The equilibrium arrival distribution contains no atoms. Indeed, if there were an atom at some time \(t\), a customer assigned to arrive at \(t\) could instead deviate to \(t-\delta\), where \(\delta>0\) is arbitrarily small, thereby obtaining a better position in the queue while changing the waiting cost only negligibly.
Moreover, there is no gap immediately after time \(0\). 
Finally, by the same reasoning as in the case without early arrivals, there is a gap immediately after each scheduled arrival time \(\mathcal{T}(k)>0\).}



\section{Algorithms and examples}\label{section:algex}

\js{From Section~\ref{section:Nash}, we obtain the following results regarding the structure of the equilibrium arrival distribution:
\begin{itemize}
    \item When early arrivals are not permitted, an atom occurs at time $0$ if no scheduled customer arrives at that time. If a scheduled customer is already scheduled to arrive at time $0$, an atom may or may not occur there.
    \item When early arrivals are permitted, the equilibrium distribution has no atoms, and the arrival density before time $0$ is nonconstant.
\end{itemize}
We next use numerical examples to illustrate these results.}
\jiesen{In our numerical experiments, we fix the service horizon at $T=10$, the service rate at $\mu=1$, and the arrival rate at $\lambda=4$. We first consider the case without early arrivals, distinguishing between whether the first scheduled customer arrives at time $0$ or after time $0$. We begin with the latter case, as it provides the basic computational framework that is subsequently adapted to the remaining settings. We then extend the analysis to allow early arrivals. For each setting, we present numerical examples illustrating the equilibrium behaviour.}

\subsection{Without early arrivals}

We distinguish between two cases according to whether the first scheduled customer arrives at time $0$.

\textbf{No scheduled customer at time $0$ ($0 \notin \mathcal{T}_s$).}
In this case, the equilibrium arrival distribution contains an atom at time $0$. At the Nash equilibrium, the expected waiting time $E_w(t)$ is constant over the support of $F_e$, and the arrival distribution satisfies $F_e(T)=1$. The equilibrium is computed as follows:
\begin{itemize}
\item Start with an initial guess for the atom size $p$ and compute the expected waiting time $E_w(0)$ for an arrival at time $0$.

\item Set $f(0)=0$ and choose a time resolution $\delta$. Compute $P_n(\delta)$ and proceed iteratively. At each time step, evaluate whether $E_w(t)\le E_w(0)$ for $t>0$. If $E_w(t)>E_w(0)$, set $f(t)=0$; otherwise, determine $f(t)$ according to~\eqref{eq:density}. Simultaneously, update the state probabilities $P_n(t)$ recursively using $f(t-\delta)$ and $P_n(t-\delta)$.

\item Compare the resulting value of $F_e(T)$ with $1$ and adjust the atom size $p$ accordingly.

\item Repeat the above steps until $F_e(T)=1$.
\end{itemize}
We can derive the validity of this computation and that the result is the unique (mixed) Nash equilibrium distribution.

\begin{proposition} \label{Uniqueness}
The value $F_e(T)$ increases monotonically as a larger atom size $p$ is chosen. Thus, there exists a unique value of $p$ such that $F_e(T)=1$. 
\end{proposition}
\begin{proof}
    Let choices $0<p_1<p_2$ be given. We denote the expected waiting times of the two choices at time $t=0$ by $E^1_{w}(0)=E^1_w,E^2_w(0)=E^2_w$, respectively. 

    We will prove this proposition by showing that the corresponding equilibrium distributions $F^1_e,F^2_e$ satisfy $ F^1_e(t)<F^2_e(t)$ for any time $0\leq t\leq T$. We will prove this statement by contradiction.

    Note that indeed $ F^1_e(0)=p_1<p_2=F^2_e(0)$.
    Assume there exists a value of $0< t\leq T$ such that $ F^1_e(t)\geq F^2_e(t)$. By continuity of $F^1_e,F^2_e$ we can take the smallest value $0< t'\leq T$ such that $F_e^1(t')=F^2_e(t')$. Let $k$ be given such that $t'\in [\mathcal{T}(k),\mathcal{T}(k+1))$.
    
    We can conclude that for both choices of atom size, the total amount of customers that entered the system (both scheduled and walk-in customers) are equal to $k$ plus a Poisson distributed variable with parameter $F_e^1(t')$. 
     Let $N^1(t'),N^2(t')$ be random variables denoting the amount of customers in the system at time $t'$ for the two respective choices $p_1,p_2$.
     As $F^1_e(t)<F^2_e(t)$ for $0<t<t'$ the amount of time that the system has been idle up to time $t'$ is stochastically smaller for the larger choice of atom size. Thus, the amount of service is stochastically larger and $ N^2(t') \preccurlyeq N^1(t')$ in the stochastic ordering. Together with the fact that $w_{\mathcal{T},k}(n,t')$ is increasing in $n$ this yields that $E_w^1<E_w^2\leq E_w^2(t')\leq E_w^1(t')$ (see Equation~(\ref{eq:Ew})). 
     $E^1_w(t)$ is continuous inside an interval before time $t'$ (with a possible upward jump at $\mathcal{T}(k)$ if $t'=\mathcal{T}(k)$). Thus,  we note that there exists a value $\epsilon>0$ such that $E^1_w<E_w^1(t)$, and therefore $f^1_e(t)=0$, for $t'-\epsilon\leq t\leq t'$ (note that this inequality also holds when $t'=\mathcal{T}(k)$).
    This gives the contradiction of the definition of $t'$ with $F^1_e(t'-\epsilon)=F^1_e(t')=F^2_e(t')\geq F^2_e(t'-\epsilon)$ for $0<t<t'$.
\end{proof}
We remark that this result does not rely on explicit expressions of the dynamics of the system. The necessary components are the homogeneity between the two types of customers, the fact that customers can't balk, and that the expected waiting time for a new walk-in customer is increasing in the current amount of customers in the queue. 


\jiesen{Although the equilibrium is characterized by analytically derived fixed-point equations, no closed-form solution exists. We therefore follow \cite{GH83} and approximate the equilibrium density \(f_e\) numerically using a finite difference method.
Let \(\delta>0\) be such that \([t,t+\delta]\subseteq [\mathcal{T}(k),\mathcal{T}(k+1))\). The system dynamics can then be approximated by
\begin{align}\label{eq:PnUpdate}
P_0(t+\delta) &= P_0(t) - \lambda f_e(t)\delta P_0(t) + \mu \delta P_1(t) + o(\delta),\\
P_n(t+\delta) &= P_n(t) - (\lambda f_e(t)+\mu)\delta P_n(t) + \lambda f_e(t)\delta P_{n-1}(t) + \mu \delta P_{n+1}(t) + o(\delta), \quad n \ge 1. \nonumber
\end{align}
All probabilities are continuous on \((\mathcal{T}(k),\mathcal{T}(k+1))\) for any \(0 \le k \le M\).}
\camiel{Moreover, we choose a truncation level $K$, defined by
\[
K = \min \left\{ k : \sum_{i=0}^k \frac{\lambda^i e^{-\lambda}}{i!} > c \right\},
\]
where $c$ is a prescribed threshold. We terminate the algorithm once the computed value of $p$ is within $\epsilon$ of the solution satisfying $F_e=1$.
Increasing $c$ and decreasing the time step $\delta$ and the tolerance $\epsilon$ improve the accuracy of the approximation. As $c\to1$ and $\delta,\epsilon\to0$, the numerical approximation converges to the true equilibrium \cite{Haviv13}. In our numerical experiments, we set $c=0.999$, $\delta=0.01$, and $\epsilon=0.001$.
}

The computation of the equilibrium distribution is summarized in Algorithm~\ref{alg1} in Appendix~\ref{app:algs}.
\camiel{
The complexity of this algorithm depends on the exact implementation, and the implementation used in our numerical experiments is available on \texttt{GitHub}. A discussion of the computational complexity is discussed after the algorithm.} 

\jiesen{

\begin{figure}[H]
\centering
\includegraphics[width=0.49\linewidth]{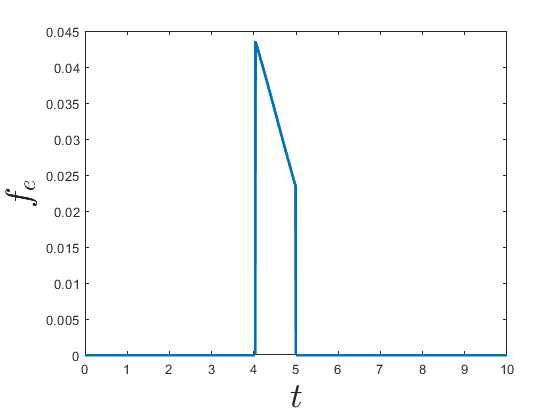}
\caption{The arrival distribution density $f_e$ when $\mathcal{T}_s=(5,6,7,8,9,10)$, in this case $p_e = 0.9678$. The other parameters are set as $\mu = 1$, $\lambda = 4$, and $T = 10$.}
\label{fig:1}
\end{figure}

We next introduce a back-loaded schedule to illustrate the equilibrium arrival behaviour under a setting where scheduled customers arrive only at later times. Specifically, consider \(\mathcal{T}_s = (5,6,7,8,9,10)\). This example is chosen to show how the absence of early scheduled arrivals affects the strategic timing decisions of walk-in customers.

The Nash equilibrium in this case yields \(p_e = 0.9678\), and the corresponding equilibrium density \(f_e\) is shown in Figure~\ref{fig:1}.
Given only the scheduled customers, the system is essentially empty during the initial time interval before the arrival of the first scheduled customer. As a result, walk-in customers have a strong incentive to concentrate their arrivals at time \(0\), leading to a large mass at the origin.
On the other hand, scheduled customers arrive in a back-loaded manner and have absolute priority, so the system becomes congested after time \(t=5\). For walk-in customers who do not arrive at time \(0\), this congestion effect makes it optimal for them to delay their arrival after time 0 but before the first scheduled customer arrives.
The resulting equilibrium pattern exhibits a mass at time \(0\), followed by a gap, and then a positive density before \(5\).}

\jiesen{{\bf A scheduled customer at time $0$ ($0 \in \mathcal{T}_s$):}
In this case, an atom may exist at $t = 0$. Let $t_0$ denote the starting point of the arrival density $f(t)$, i.e.,
\begin{equation}\label{def:t0}
t_0 = \inf\{t \mid f(t) > 0\}.
\end{equation}
The function $f(t)$ for $t \in [t_0, T]$ can then be computed iteratively using~\eqref{eq:density}. In the absence of an atom, a larger $t_0$ leads to a smaller $E_w(t_0)$ and consequently a smaller $f(t)$ over the interval $[t_0, T]$. Thus, increasing $t_0$ results in a smaller $F_e(T)$.
If $F(T) < 1$ even when $t_0 = 0$, this indicates the presence of an atom at $t = 0$. \camiel{We note that $F(T)$ decreases monotonically as $t_0$ increases, by the same argument as in Proposition~\ref{Uniqueness}. This monotonicity also confirms the uniqueness of the Nash equilibrium distribution.} 

Based on this analysis, the Nash equilibrium is computed as follows:
\begin{itemize}
\item First, we compute $F(T)$ under the assumption that $t_0=0$. If $F(T)<1$, we apply Algorithm~\ref{alg1} to compute $p_e$ and $f_e(t)$.
\item Otherwise, we initialize $t_0=\mathcal{T}(2)$ and compute the corresponding arrival density $f(t)$. We then compare the resulting value of $F(T)$ with $1$.
\item If $F(T)<1$, we increase $t_0$ to the next scheduled arrival time and repeat the previous step. Otherwise, update $t_0$ using a bisection search based on the value of $F(T)$, and repeat the previous step until $F(T)=1$.
\end{itemize}
This procedure is summarized in Algorithm~\ref{alg2} in Appendix~\ref{app:algs}. A discussion of its computational complexity is provided after the algorithm.}

\jiesen{In our numerical experiments, we consider two appointment schedules, each consisting of six scheduled customers. Both schedules include an initial appointment at time \(0\). We vary the spacing of subsequent appointments in order to study how different scheduling patterns affect equilibrium behaviour.
Specifically, we consider \js{
\[
\mathcal{T}_s^{(1)}=(0,1,2,3,4,5)
\quad\text{and}\quad
\mathcal{T}_s^{(2)}=(0,2,4,6,8,10).
\]
}
The first schedule places appointments relatively densely at the beginning of the planning horizon, whereas the second is evenly spaced over time. This contrast allows us to compare a front-loaded schedule with a uniformly distributed one.
The corresponding equilibrium arrival densities are shown in Figure~\ref{fig:2}.

\begin{figure}[H]
\centering

\subcaptionbox{\js{$\mathcal{T}_s^{(1)} = (0, 1, 2, 3, 4, 5):\; p_e = 0 \,$}}
{\includegraphics[width=0.49\linewidth]{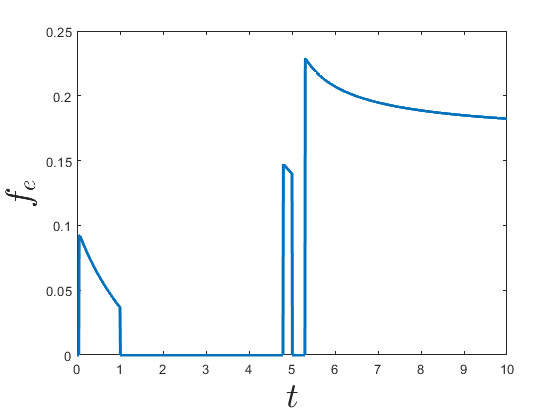}}
\subcaptionbox{$\mathcal{T}_s^{(2)} = (0, 2, 4, 6, 8, 10):\; p_e = 0.3394$}
{\includegraphics[width=0.49\linewidth]{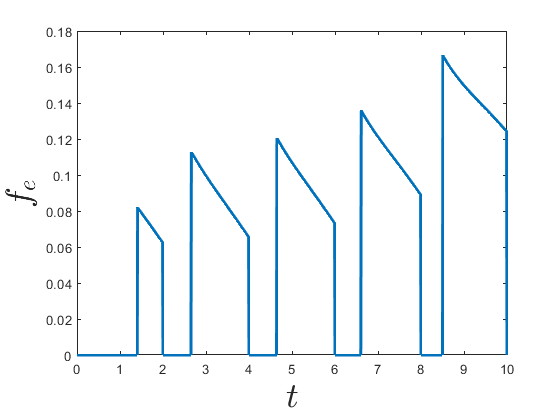}}

\caption{The arrival distribution density $f_e$. The other parameters are set as $\mu = 1$, $\lambda = 4$, and $T = 10$.}
\label{fig:2}
\end{figure}
It can be seen that, in the front-loaded case, most walk-in customers choose to arrive after all scheduled customers have arrived, although some still arrive in the gaps between scheduled arrivals. Because there is a scheduled customer at time $0$ with absolute priority, arriving exactly at time $0$ is no longer an optimal strategy. Consequently, the equilibrium contains no atom at time $0$ ($p_e=0$), and the arrival density starts at $t_0=0.05\approx0$. Customers therefore prefer to arrive shortly after time $0$ rather than exactly at time $0$. By delaying their arrival slightly, customers retain the same queue position but with a positive probability that the scheduled customer has already departed. It is also worth noting that, under our model, a walk-in customer can be delayed by at most six scheduled customers. This worst-case scenario occurs when every scheduled customer arrives before the previous scheduled customer has completed service. Consequently, the expected waiting time of a walk-in customer is bounded from above by the service time of six scheduled customers together with the expected service time of all walk-in customers that arrived before.

For the evenly spaced schedule, it is clear that a gap exists after each scheduled arrival. In the case without early arrivals and with the first scheduled arrival at time zero, the structure is more delicate and depends on whether an atom is present at time zero. If no atom is present, there may be a gap immediately after time zero.
}

\clearpage

\subsection{Including early arrivals}

\jiesen{In this case, it is possible that no customer has an incentive to arrive earlier, even when early arrival is allowed. Specifically, if the equilibrium arrival distribution of walk-in customers in the model without early arrivals has no atom at time zero, it also remains an equilibrium arrival distribution when early arrivals are allowed.
This follows from the fact that \(E_w(0) \geq E_w\), which implies \(E_w(t) > E_w\) for all \(t < 0\). Moreover, when early arrivals are allowed, overlaps may arise in later parts of the horizon, as illustrated in the numerical example below.

We take \(t_0\) as defined in a manner similar to Equation~\eqref{def:t0}, with the interpretation that it serves as the initial point of the density. The equilibrium is computed according to the following procedure:

\begin{itemize}
    \item Initialize with \(t_0=0\) and compute \(F(T)\). If \(F(T)<1\), decrease \(t_0\) to a negative value.

    \item Check whether \(E_w(t)\le E_w(t_0)\) for \(t>t_0\). If \(E_w(t)>E_w(t_0)\), set \(f(t)=0\); otherwise, update \(f(t)\) according to \eqref{eq:earlyfe} for \(t<0\) and according to \eqref{eq:density} for \(t\ge 0\).

    \item Compute the resulting equilibrium mass \(F_e\) and adjust \(t_0\) accordingly.

    \item Repeat the above steps until \(F_e=1\).
\end{itemize}

The corresponding algorithm is obtained by modifying Algorithm~\ref{alg2} to incorporate the steps for determining \(t_0\) when it is negative. The complete procedure is presented as Algorithm~\ref{alg3} in Appendix~\ref{app:algs}. As in Proposition~\ref{Uniqueness}, \(F(T)\) is monotone in \(t_0\), implying uniqueness of the Nash equilibrium distribution.
}

\begin{figure}[H]
\centering
\subcaptionbox{$\mathcal{T}^{(1)}_s = (0, 2, 4, 6, 8, 10):\; p_e = 0.3394\,$}
{\includegraphics[width=0.49\linewidth]{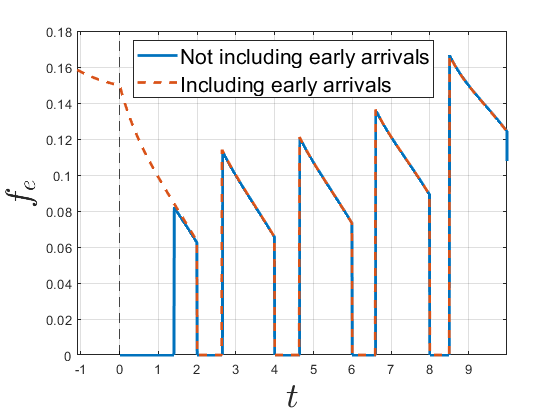}}
\subcaptionbox{$\mathcal{T}^{(2)}_s = (1, 2, 4, 6, 8, 10):\; p_e = 0.5557$}
{\includegraphics[width=0.49\linewidth]{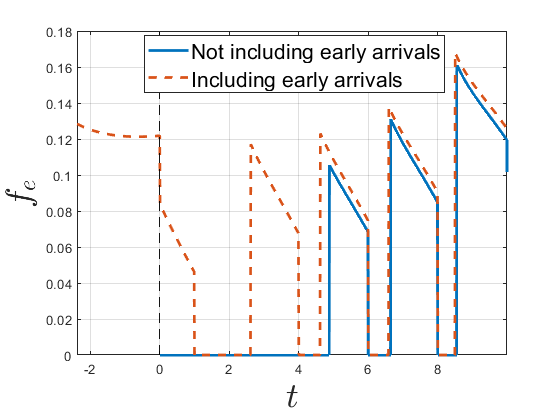}}
\caption{Equilibrium arrival density \(f_e\) with and without early arrivals for two appointment schedules: \(\mathcal{T}_s^{(1)}=(0,2,4,6,8,10)\), for which \(\mathcal{T}_s^{(1)}(1)=0\), and \(\mathcal{T}_s^{(2)}=\js{(1, 2, 4, 6, 8, 10)}\), for which \(\mathcal{T}_s^{(2)}(1)>0\). The parameters are \(\mu=1\), \(\lambda=4\), and \(T=10\).}
\label{fig:3}
\end{figure}

\jiesen{In the numerical analysis with early arrivals, we consider two appointment schedules,
\[
\mathcal{T}_s^{(1)}=(0,2,4,6,8,10),\qquad
\mathcal{T}_s^{(2)}=(1,2,4,6,8,10).
\]
We investigate how the inclusion of a time-zero appointment affects the equilibrium early arrival distribution.
For each schedule, we plot the corresponding equilibrium density \(f_e\) together with the equilibrium density obtained in the case without early arrivals; see Figures~\ref{fig:3}(a) and~\ref{fig:3}(b).

When \js{the first scheduled customer arrives at time zero}, 
the equilibrium density \(f_e\) is continuous at time \(0\), although its derivative is discontinuous. Moreover, the continuous parts of the equilibrium distributions with and without early arrivals have a near exact overlap in the later part of the planning horizon, and the corresponding expected waiting times in this case can be nearly indistinguishable.

In contrast, when \js{\js{the first scheduled customer arrives after time zero}}, 
the equilibrium density \(f_e\) exhibits a discontinuity at time \(0\). Furthermore, the continuous equilibrium density under early arrivals is strictly higher than its counterpart without early arrivals over most of the planning horizon.
 This difference can also be explained as the early difference between the models is not damped by the first scheduled customer being in service. Furthermore, in the model with early arrivals a part of the waiting time of walk-in customers can be the waiting time for the system to commence. 
As a result, allowing early arrivals leads to a larger equilibrium expected waiting time. We observed the same qualitative behaviour in several additional examples.}

\jiesen{
\section{Scheduling under behavioural uncertainty} \label{section:schedule}

In Section \ref{section:Nash}, we derived the equilibrium arrival distribution. A natural follow-up question is how the scheduling policy influences this distribution, and what the optimal schedule is when both idle time and waiting time are taken into account.

The work \cite{Wang93} addresses the problem of scheduling a finite number of customer arrivals to a single-server system in which customers arrive by appointment only, with the objective of minimizing a weighted combination of customer delay and server completion time. Our work considers a similar objective; however, a key distinction is that we incorporate strategically behaving walk-in customers who choose their arrival times endogenously. These strategic arrivals, in turn, affect the waiting times of scheduled customers and thus alter the overall system performance.

In practice, most clinics or medical centers adopt evenly spaced appointment schedules, which is a rather coarse approach. We will later show numerically that although this policy is not optimal, it does not lead to a large performance loss. 
In reality, both the arrival rate of walk-ins and the fraction of strategic versus non-strategic patients are typically not directly observable. Thus, it is important that the chosen schedule is robust with respect to actual customer behaviour.
We observe that the equally spaced scheduling policy is a reasonably robust and practically efficient heuristic: though not optimal, it performs reasonably well under such uncertainty.

\js{Specifically, we consider a cost function with three components: the expected waiting time of walk-in customers, the expected waiting time of scheduled customers, and the server's expected idle time. Minimizing the server's idle time is equivalent to minimizing the expected workload remaining when the system closes, corresponding to the makespan criterion in this setting.
Let $E_I(\mathcal{T})$ denote the server's expected idle time, $E_w(\mathcal{T})$ the expected waiting time of an individual walk-in customer, and $E_s(\mathcal{T})$ the total expected waiting time of all scheduled customers. Since walk-in customers have the same expected waiting time under a symmetric strategy, the total expected waiting time of walk-in customers is $\lambda E_w(\mathcal{T})$.
The objective function under strategic behaviour is defined as
\[
\Phi(\mathcal{T})
:= (1-\gamma)\left(
E_s(\mathcal{T})+\lambda E_w(\mathcal{T})
\right)
+\gamma E_I(\mathcal{T}),
\]
}
where $\gamma\in[0,1]$ captures the relative importance of reducing server idle time versus minimizing customer waiting times.

When customers are naive, their arrival distribution is assumed to be uniform, whereas when they are strategic, their behaviour corresponds to the Nash equilibrium arrival distribution derived in Section~\ref{section:Nash}.
We assume that $M$ is fixed, i.e., the number of scheduled customers is given, and we numerically analyze three types of scheduling policies:
\begin{itemize}
\item Case I: evenly spaced appointment schedules;
\item Case II: the optimal schedule under the assumption that the walk-in rate $\lambda$ is known and customers behave naively;
\item Case III: the optimal schedule under the assumption that $\lambda$ is known and customers behave strategically.
\end{itemize}

Our goal is to quantify the performance loss of Cases I and II relative to Case III  when customers behave strategically. Importantly, all performance metrics are evaluated under strategic customer behaviour; the distinction between the cases lies only in the assumptions used by the service provider when designing the schedule.

\textbf{Structure of the optimal schedule without walk-ins}. 
When there are no walk-in customers, placing greater emphasis on minimizing idle time leads to a front-loaded optimal strategy. In the extreme case $\gamma = 1$, the objective reduces to minimizing idle time alone, and the optimal policy concentrates all appointments at the beginning of the horizon, thereby maximizing the likelihood that the system remains continuously busy until all customers are served. As $\gamma$ decreases, the optimal schedule becomes more dispersed over time and, when $\gamma = 0$, exhibits a \emph{dome-shaped} pattern \cite{Wang93} , where inter-arrival times are larger at the two ends and smaller in the middle. For example, when $\lambda = 0$, $M = 6$, and $T = 10$, the optimal schedule for $\gamma = 0$ is given by
\[
\big[0, \quad 1.87, \quad 4.05, \quad 6.24, \quad 8.37, \quad 10 \big].
\]
To minimise waiting time, the first customer is always scheduled at the start of the interval, $t_1 = 0$, and the last customer at the end, $t_M = T$.
This dome-shaped structure reflects non-uniform spacing of appointments: inter-arrival times are smallest at the beginning and end, and largest in the middle. In the early phase, gaps are small since the system is initially empty and congestion has not yet built up. In the middle phase, gaps widen to mitigate the accumulation of delays, which tend to peak in this region. In the late phase, the final gap is smaller, since delays for the last customer do not propagate further; it is therefore less costly to allocate more waiting time to the last arrival in exchange for reducing waiting times earlier in the schedule. 

\textbf{Performance comparison of scheduling policies}. Given a schedule $\mathcal{T}_s$, we can compute the expected cost $\Phi(\mathcal{T}_s)$. From the service provider’s perspective, neither the arrival rate of walk-ins $\lambda$, nor the extent to which customers behave strategically is directly observable.
To assess the impact of this information gap, we consider two benchmark policies: the evenly spaced schedule (Case I), and the schedule obtained under the assumption that walk-in arrivals are uniformly distributed (Case II).
In contrast, we also derive the optimal schedule and the corresponding cost under the assumption that customers behave strategically (Case III). Comparing these cases quantifies the value of accounting for strategic behaviour and having accurate knowledge of the system parameters.
In many healthcare settings, reducing idle time is a primary operational concern. Although $\gamma$ may take different values, we focus on $\gamma = 0.9$ in the following analysis.

The optimisation problem is challenging due to the presence of multiple local optima and the indirect relationship between the schedule and the objective function. In addition, the use of a discretisation parameter $\delta$ in the approximation further complicates the optimisation landscape.

We solve the problem using MATLAB’s \texttt{patternsearch} algorithm, a derivative-free method based on direct search. Starting from an initial guess, the algorithm iteratively explores a set of candidate points in a predefined pattern around the current solution and updates the solution whenever an improvement is found. The search proceeds with adaptive step sizes until convergence criteria are satisfied. This approach is well suited for problems where the objective function may be non-smooth or lacks reliable gradient information.

The computation time for Case III for $T=10$ and $M=6$ is substantial and increases with the value of $\lambda$. For example, each run takes approximately 30 hours when $\lambda = 2$, and about 40 hours when $\lambda = 6$. For each value of $\lambda$, we test five different initial values for the scheduling decision and the associated optimization procedure.

For $\lambda = 2$ and $\lambda = 4$, all runs converge to the same solution. However, for $\lambda = 6$, the algorithm does not consistently converge to a single value; instead, the obtained solutions lie in the range $[3.7168, 3.8197]$. Since our computation is based on an approximation scheme, achieving higher accuracy would require reducing $\delta$ and $\varepsilon$, which would further increase the computational time.
Therefore, for $\lambda = 6$, we report $3.7168$ as the optimal value, although we acknowledge that it may not be the global optimum.

\begin{table}[H]
    \centering
    \begin{tabular}{||l|l|l|l|c||}
    \hline
    $\lambda$ & Assumption & Optimal schedule & Cost & $\frac{(\Phi(\mathcal{T})-\Phi(\mathcal{T}^{*(3)})}{\Phi(\mathcal{T}^{*(3)})}$ \\
    \hline
    \hline
    
    \multirow{3}{*}{$2$}
    & evenly spaced & $\mathcal{T}_s^{(1)}$ & $ 4.0214$ & 16.11\%\\
    & naive & $\mathcal{T}_s^{*(2)} = [0, \, 0.51, \, 1.85, \, 3.49, \, 5.25, \, 7.09]$ & $3.5288$ & 1.89\%\\
    & strategic & $\mathcal{T}_s^{*(3)} = [0.04, \,  1.22, \,   2.49, \,  3.75, \,  5.16, \,   6.32]$ & $3.4634$ & 0 \\
    \hline
    
    \multirow{3}{*}{$4$}
    & evenly spaced & $\mathcal{T}_s^{(1)}$ & $3.7013$ & 11.72\%\\
    & naive & $\mathcal{T}_s^{*(2)} = [0, \, 0.31, \, 1.57, \, 3.38, \, 5.42, \, 7.49]$ & $3.3376$ & 0.74\% \\
    & strategic & $\mathcal{T}_s^{*(3)} = [0.08, \, 1.29, \, 2.47, \, 3.70, \, 4.65, \, 5.96]$ & $3.3131$  & 0 \\
    \hline
    
    \multirow{3}{*}{$6$}
    & evenly spaced & $\mathcal{T}_s^{(1)}$ & $4.1437$ & 11.49\%\\
    & naive & $\mathcal{T}_s^{*(2)} = [0, \, 0.23, \, 1.45, \, 3.64, \, 6.23, \, 8.47]$ & $3.8104$ & 2.52\% \\
    & strategic & $\mathcal{T}_s^{*(3)} = [0, \, 0.15, \, 0.97, \, 2.13, \, 3.60, \, 5.49]$ & $3.7168$ & 0 \\
    \hline
    \end{tabular}
    \caption{Comparison of scheduling policies for $M=6$, $T=10$, and $\gamma=0.9$. The schedule 
$\mathcal{T}_s^{(1)} = [0,\,2,\,4,\,6,\,8,\,10]$ is the evenly spaced schedule, while 
$\mathcal{T}_s^{*(2)}$ and $\mathcal{T}_s^{*(3)}$ denote the optimal schedules obtained under the assumptions of uniform (naive) arrivals and strategic arrivals, respectively. 
We assume that customers are in fact strategic and evaluate the performance loss resulting from using a schedule optimized under the incorrect assumption of uniform arrivals. 
The corresponding costs are reported in the order 
$\Phi(\mathcal{T}_s^{(1)}), \, \Phi(\mathcal{T}_s^{*(2)}), \, \Phi(\mathcal{T}_s^{*(3)})$.}
    \label{tab:cost}
\end{table}


Across the three cases, the evenly spaced schedule yields the highest cost, but the gap is reasonably small. Given its simplicity, this helps explain why it is widely adopted in practice. In contrast, the schedules derived under uniform and strategic assumptions (Cases II and III) exhibit very similar performance, indicating that the optimisation is robust to misspecification of customer behaviour.

An additional observation is that, under evenly spaced schedules, strategic behaviour leads to higher waiting times for walk-in customers, implying that strategic behaviour can be detrimental from the walk-in customers’ perspective in this setting. Moreover, although Cases II and III yield nearly identical total costs, the expected waiting time of walk-in customers is lower in Case III. This observation suggests that strategic customers may benefit when the planner explicitly recognizes and incorporates their equilibrium behaviour into the scheduling design. In this sense, strategic walk-in customers may have an incentive to make their strategic nature known, as the resulting schedule is better aligned with their preferences.

\js{Note that the total cost is higher when $\lambda=4$ than when $\lambda=2$. This is because our objective function assigns a higher weight to idle time. When $\lambda=2$, the total service capacity exceeds demand, resulting in more idle time and hence a higher total cost.}
}

\section{Conclusion}
\label{section:conclusion}

This paper investigates the dynamics of a single-server queueing system that accommodates both scheduled and strategic walk-in customers, operating under a first-come, first-served (FCFS) discipline with non-preemptive priority for scheduled customers. Scheduled customers induce an exogenous, time-varying congestion profile that reshapes the strategic arrival behaviour of walk-in customers.
By modeling the strategic behaviour of walk-in customers, who aim to minimise their expected waiting time, we derive the Nash equilibrium arrival distribution and analyse how different appointment schedules impact system performance. Our findings reveal the complex interplay between scheduled and walk-in customers, particularly how the timing of scheduled appointments influences walk-in arrival patterns and overall system efficiency.
Our main contributions are as follows:
\begin{itemize}\item[(i)] We characterise the equilibrium arrival behaviour of strategic walk-in customers in a queueing system with exogenous scheduled demand.
\item[(ii)] We quantify the impact of appointment scheduling on equilibrium outcomes and system-level performance trade-offs.
\end{itemize}

Our numerical results show that scheduling appointments earlier, reduces waiting times for scheduled customers, but increases waiting times for walk-in customers, whereas scheduling appointments later, benefits walk-in customers at the expense of higher server idle time. Furthermore, while aggregate performance measures under the equilibrium arrival pattern may differ only moderately from those obtained under a uniform random-arrival assumption, strategic behaviour substantially alters the temporal distribution of walk-in arrivals. 
\camiel{Finally, we conclude that the reasonable performance loss by using evenly spaced schedules with the presence of strategic customers leads to a debatable trade-off between performance and practicality.}


\section{Acknowledgment}
The authors would like to thank Liron Ravner for his valuable comments and advice. This research was supported by the European Union’s Horizon 2020 research and innovation programme under the Marie Sklodowska-Curie grant agreement no.\ 945045, and by the NWO Gravitation project NETWORKS under grant agreement no.\ 024.002.003. \includegraphics[height=1em]{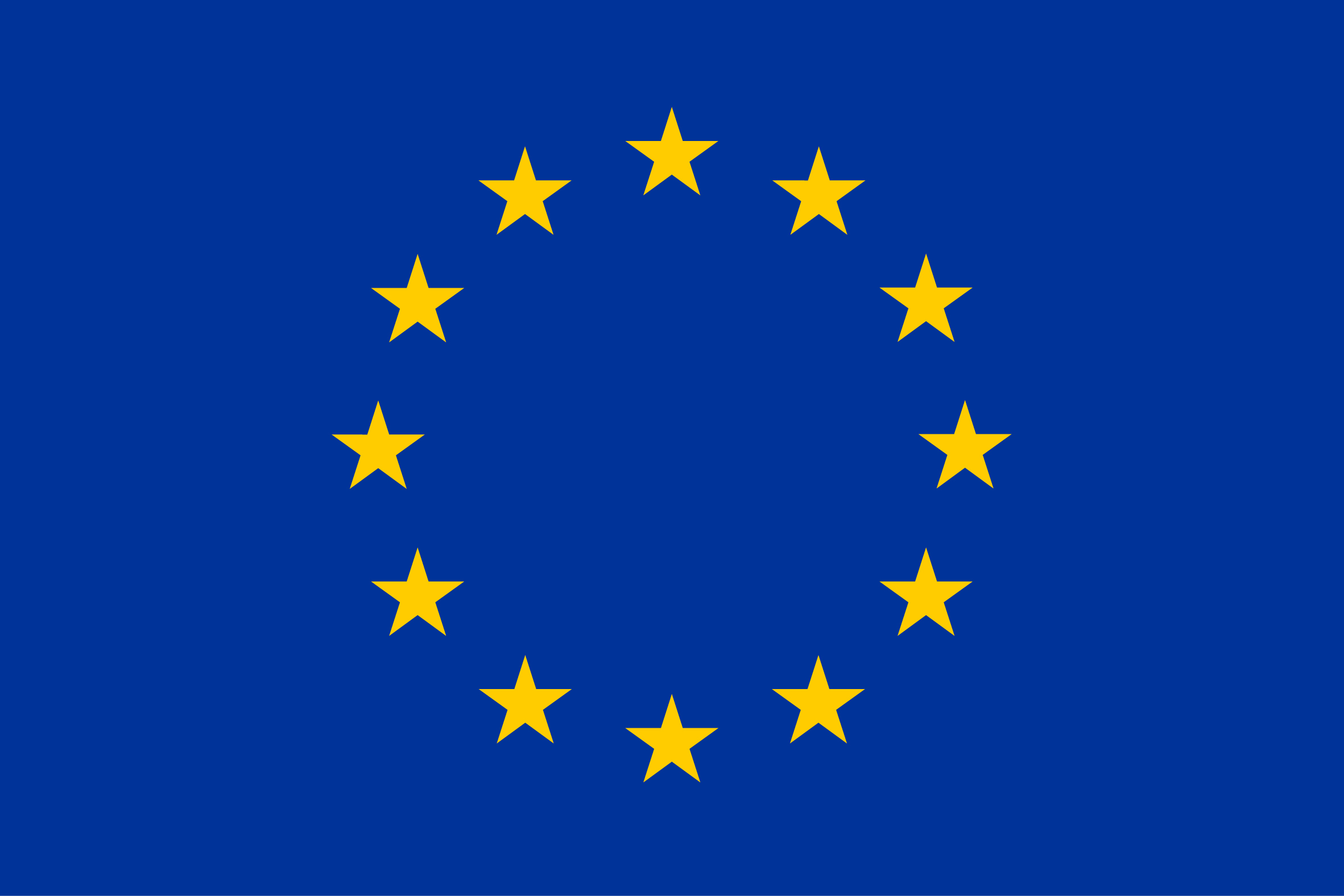}

\clearpage

\appendix
\section{\js{Derivative of $\sum_{n=0}^{\infty} g(n,t)P_n(t)$ with respect to $t$}}
\label{Derivatives}

\js{For a function $g:\mathbb{N}_0\to\mathbb{R}$, and for
$t\neq \mathcal{T}(i)$, $i=1,2,\ldots,M$, the derivative of
$\sum_{n=0}^{\infty}g(n)P_n(t)$ with respect to $t > 0$ is given by
\begin{align*}
    \frac{\partial}{\partial t}\sum_{n=0}^{\infty}g(n)P_n(t)
    &= \lim_{\delta\to 0}
    \frac{\sum_{n=0}^{\infty}g(n)P_n(t+\delta)
    -\sum_{n=0}^{\infty}g(n)P_n(t)}{\delta}\\
    &= \sum_{n=0}^{\infty}
    \bigl(g(n+1)-g(n)\bigr)
    \bigl(P_n(t)\lambda f(t)-P_{n+1}(t)\mu\bigr).
\end{align*}

We now extend this result to the case in which $g$ also depends on
time, that is, $g:\mathbb{N}_0\times\mathbb{R}\to\mathbb{R}$.
Suppose that, for each $n\in\mathbb{N}_0$, $g(n,t)$ is differentiable
with respect to $t > 0$. Then, for $t\neq \mathcal{T}(i)$,
$i=1,2,\ldots,M$, the derivative of
$\sum_{n=0}^{\infty}g(n,t)P_n(t)$ is given by
\begin{align}
\frac{\partial}{\partial t}\sum_{n=0}^{\infty}g(n,t)P_n(t) = &\lim_{\delta\to 0}
\frac{\sum_{n=0}^{\infty}g(n,t+\delta)P_n(t+\delta)
-\sum_{n=0}^{\infty}g(n,t)P_n(t)}{\delta}
\nonumber\\
=&\lim_{\delta\to 0}\Bigg[
\sum_{n=0}^{\infty}P_n(t)
\frac{g(n,t+\delta)-g(n,t)}{\delta}
\nonumber\\
&\qquad\quad+
\lambda f(t)\sum_{n=0}^{\infty}P_n(t)
\bigl(g(n+1,t+\delta)-g(n,t+\delta)\bigr)
\nonumber\\
&\qquad\quad-
\mu\sum_{n=0}^{\infty}P_{n+1}(t)
\bigl(g(n+1,t+\delta)-g(n,t+\delta)\bigr)
\Bigg]
\nonumber\\
&=\sum_{n=0}^{\infty}P_n(t)
\frac{\partial g(n,t)}{\partial t}
+\lambda f(t)\sum_{n=0}^{\infty}P_n(t)
\bigl(g(n+1,t)-g(n,t)\bigr)
\nonumber\\
&\qquad
-\mu\sum_{n=0}^{\infty}P_{n+1}(t)
\bigl(g(n+1,t)-g(n,t)\bigr).
\label{eq:DeGP}
\end{align}

When $t<0$, the system has not yet commenced service, so no departures occur. Hence, the derivative of $\sum_{n=0}^{\infty}g(n)P_n(t)$ is given by
\begin{equation}\label{eq:deNPEarly1}
    \sum_{n=0}^{\infty}
    \bigl(g(n+1)-g(n)\bigr)P_n(t)\lambda f(t).
\end{equation}
Similarly, the derivative of $\sum_{n=0}^{\infty}g(n,t)P_n(t)$ is
\begin{equation}\label{eq:deGPEarly2}
    \sum_{n=0}^{\infty}P_n(t)\frac{\partial g(n,t)}{\partial t}
    +\lambda f(t)\sum_{n=0}^{\infty}P_n(t)
    \bigl(g(n+1,t)-g(n,t)\bigr).
\end{equation}
}

\newpage

\section{Algorithms}
\label{app:algs}

\begin{algorithm}[H]
\caption{Computation of the Nash equilibrium when $0\notin\mathcal{T}_s$}
\label{alg1}
\begin{algorithmic}[1]
\State \textbf{Input:} $\lambda,\mu,T,\mathcal{T},\delta,c,\epsilon$
\State \textbf{Initialize:} $p_L\gets0,\; p_U\gets1$
\While{$p_U-p_L>\epsilon$}
    \State $p\gets(p_L+p_U)/2$
    \State Initialize $P_n(0)$ using the atom size $p$
    \State Compute $E_w(0)$
    \State Set $f(0)\gets0$
    \For{$k=0$ to $M$}
        \For{$t=\mathcal{T}(k)$ to $\mathcal{T}(k+1)$}
            \State Update $P_n(t)$ and compute $E_w(t)$
            \If{$E_w(t)>E_w(0)$}
                \State Set $f(t)\gets0$
            \Else
                \State Compute $f(t)$ using~\eqref{eq:density}
            \EndIf
        \EndFor
    \EndFor
    \State Compute the corresponding value of $F(T)$
    \If{$F(T)>1$}
        \State $p_U\gets p$
    \Else
        \State $p_L\gets p$
    \EndIf
\EndWhile
\State \textbf{Output:} $(p_e,f_e)$
\end{algorithmic}
\end{algorithm}

\camiel{In order to save computation time, we first compute a matrix with all values $\Tilde{w}_{\mathcal{T},k}(n)$ for $0\leq k \leq M$ and $0\leq n\leq K$. This can be done by starting with $k=M$ and computing the values $\Tilde{w}_{\mathcal{T},k}(n)$ for $0\leq n\leq K$ with a computation time of $\mathcal{O}(k)$ per element (cf. Equation~(\ref{eq:wkt})). Iterating backwards in $k$ results in a total computation time of $\mathcal{O}(M\cdot K^2)$ for the whole matrix. We note that, in order to calculate the whole matrix correctly, it is necessary to include calculation of values $\Tilde{w}_{\mathcal{T},k}(n)$ for $0\leq n\leq K+k$. This does not increase the order of the computation time.
The binary search of $p$ will terminate after $\mathcal{O}(\log(1/\epsilon))$ iterations.

Each iteration starts with the
calculation of $E_w(0)$ (cf. Equation~(\ref{eq:Ew0})) of $\mathcal{O}(K^2)$ time and proceeds with $\mathcal{O}(T/\delta)$ time steps. Each time step consists of an update of values $P_n(t)$ taking $\mathcal{O}(K)$ time (cf. Equation~(\ref{eq:PnUpdate}) and Equation~(\ref{eq:PschedArr}) in case of a scheduled customer), calculation of $E_w(t)$ taking $\mathcal{O}(K^2)$ time (cf. Equations~(\ref{eq:wkt}) and (\ref{eq:Ewt})) and finally the computation of $f(t)$ taking $\mathcal{O}(K^2)$ time (cf. Equation~(\ref{eq:density}) with Equations~(\ref{eq:wkt}) and (\ref{eq:derwkt})).

We can conclude that the total complexity of the algorithm is given by\\ $\mathcal{O}(K^2(M+ \log(1/\epsilon)\cdot T/\delta))$.}

\begin{algorithm}[H]
\caption{Computation of the Nash equilibrium when $0\in\mathcal{T}_s$}
\label{alg2}
\begin{algorithmic}[1]
\State \textbf{Input:} $\lambda,\mu,T,\mathcal{T},\delta,c,\epsilon$
\State Compute $F(T)$ assuming $t_0=0$.
\If{$F(T)<1$}
    \State Apply Algorithm~\ref{alg1} to compute $p_e$ and $f_e(t)$.
\Else
    \State Initialize $l\gets1$ and $t_0\gets\mathcal{T}(2)$.
    \State Compute $f(t)$ and the corresponding value of $F(T)$.
    \While{$F(T)>1$}
        \State $l\gets l+1$.
        \State $t_0\gets\mathcal{T}(l)$.
        \State Compute $f(t)$ and the corresponding value of $F(T)$.
    \EndWhile
    \State Set $t_L\gets\mathcal{T}(l-1)$ and $t_U\gets\mathcal{T}(l)$.
    \While{$|F(T)-1|>\epsilon$}
        \State $t_0\gets (t_L+t_U)/2$.
        \State Compute $f(t)$ using~\eqref{eq:density}.
        \State Compute the corresponding value of $F(T)$.
        \If{$F(T)>1$}
            \State $t_L\gets t_0$.
        \Else
            \State $t_U\gets t_0$.
        \EndIf
    \EndWhile
\EndIf
\State \textbf{Output:} $(p_e,f_e)$.
\end{algorithmic}
\end{algorithm}

\camiel{
Regarding the complexity of this algorithm we note that we either use Algorithm~\ref{alg1}
with complexity $\mathcal{O}(K^2(M+ \log(1/\epsilon)\cdot T/\delta))$, or $t_0\geq 0$ and we consider $\mathcal{O}(K+\log(T/\delta))$ guesses of $t_0$ (the $K$ times of scheduled customers and the binary search). Each guess, again results in a complexity of $\mathcal{O}(K^2\cdot T/\delta))$. Thus, the complexity of Algorithm~\ref{alg2} is given by $\mathcal{O}(K^2(M+ (\log(1/\epsilon)+K+\log(T/\delta) )\cdot T/\delta))$.}

\begin{algorithm}[H]
\caption{Computation of the Nash equilibrium with early arrivals}
\label{alg3}
\begin{algorithmic}[1]
\State \textbf{Input:} $\lambda,\mu,T,\mathcal{T},\delta,c,\epsilon$
\State Compute $F(T)$ with $t_0\gets 0$.
\If{$F(T)>1$}
    \State Apply Algorithm~\ref{alg2} to compute the equilibrium distribution.
\Else
    \State \textbf{Initialize:} Find a negative initial value \(t_0\) by decrementing it by 1 at each step until \(F(t_0)\ge 1\). Then set \(t_L \leftarrow t_0\) and \(t_U \leftarrow 0\).
    \While{$t_U-t_L>\delta$ \textbf{or} $|F(T)-1|>\epsilon$}
        \State $P_n(t_0)\gets[1,0,0,\ldots]$.
        \For{$t=t_0$ \textbf{to} $0$}
            \State Compute $f(t)$ using \eqref{eq:earlyfe}.
            \State Update $P_n(t)$.
        \EndFor
        \For{$k=0$ \textbf{to} $M$}
            \For{$t=\mathcal{T}(k)$ \textbf{to} $\mathcal{T}(k+1)$}
                \State Compute $E_w(t)$.
                \If{$E_w(t)\le E_w(t_0)$}
                    \State Compute $f(t)$ using \eqref{eq:density}.
                \Else
                    \State Set $f(t)\gets 0$.
                \EndIf
                \State Update $P_n(t)$.
            \EndFor
        \EndFor
        \State Compute $F(T)$.
        \State Update $t_0$ using a bisection search.
    \EndWhile
\EndIf
\State \textbf{Output:} $f_e$.
\end{algorithmic}
\end{algorithm}

\camiel{For the complexity of this algorithm we need to find an order of the amount of guesses needed for $t_0$, for this we need to find a lower bound $L\leq t_0$. To find this lower bound we observe that $w_{\mathcal{T},\mathbbm{1}_{\mathcal{T}(1)=0}}(n+1+\mathbbm{1}_{\mathcal{T}(1)=0},0)-w_{\mathcal{T},\mathbbm{1}_{\mathcal{T}(1)=0}}(n+\mathbbm{1}_{\mathcal{T}(1)=0},0)\leq M+1$. Thus, Equation~(\ref{eq:earlyfe}) yields $f_e(t)\geq 1/(\lambda(M+1))$ for $t_0\leq t<0$. Thus, $F_e(0)\geq -t_0/(\lambda(M+1))$. Thus, $-\lambda(M+1)$ is a sufficient lower bound of $t_0$.

As a result, there are $\mathcal{O}(\log((\lambda(M+1)+T)/\delta))$ time steps taken.
In our implementation, $t_0$ is reduced by 1 at most $\lambda(M+1)$ times after which a binary search of $\mathcal{O}(1/
\delta)$ searches is conducted. 
This leads to a time complexity of 
$\mathcal{O}(K^2(M+ \log(T/\delta) \cdot T/\delta+\lambda(M+1)\log(1/\delta)(\lambda(M+1)+T)/\delta))$.
The worst case complexity can be reduced by considering a binary search on $t_0$ for $t_0\in [-\lambda(M+1),0)$ instead of repetitive reduction of $t_0$ by 1. The reasoning for the current implementation is that the lower bound $-\lambda(M+1)$ often is rather pessimistic.}

\end{document}